\documentclass[11pt]{amsart}
\usepackage{graphicx}
\usepackage{amssymb}
\newcommand{\C}{\mathcal{C}}

\newtheorem{Def}{Definition}
\newtheorem{Lem}{Lemma}

\newtheorem{Thm}{Theorem}

\newtheorem{Rem}{Remark}
\newenvironment{Pf}{ Proof.}{\(\square\)}

 \title{An observation on Asanov's Unicorn metrics}
\author{ Csaba Vincze}  
\address{ Institute of Mathematics, University of Debrecen\\ H-4010 Debrecen, P.O. Box 12 Hungary }
\email{csvincze@science.unideb.hu}
\keywords{Finsler spaces, Conformality, Finsleroid-Finsler metrics, Generalized Berwald manifolds.}
\subjclass{53C60} 
\thanks{ Supported by the University of Debrecen's internal research project RH/885/2013.}

\begin{document}

\begin{abstract}
Finsleroid-Finsler metrics form an important class of singular (y-local) Finslerian metrics. They were introduced by G. S. Asanov in 2006. As a special case Asanov produced examples of  Landsberg spaces of dimension at least three that are not of Berwald type. These are called Unicorns \cite{B}. The existence of regular (y - global) Landsberg metrics that are not of Berwald type is an open problem up to this day. In this paper we prove that Asanov's Unicorns belong to the class of generalized Berwald manifolds. More precisely we prove the following theorems: a Finsleroid-Finsler space is a generalized Berwald space if and only if the Finsleroid charge is constant. Especially a Finsleroid-Finsler space is a Landsberg space if and only if it is a generalized Berwald manifold with a semi-symmetric compatible linear connection. 
\end{abstract}

\maketitle
               
\section*{Introduction}

Finsleroid-Finsler metrics were introduced by G. S. Asanov in 2006. As a special case Asanov produced singular (y - local) examples of  Landsberg spaces of dimension at least three that are not of Berwald type. The existence of regular (y - global) Landsberg metrics that are not of Berwald type is an open problem up to this day, see D. Bao \cite{B}. The cronology of the basic steps:
\begin{itemize}
\item 1998 - the central symmetric version of the Finsleroid-Finsler metric in \cite{Asanov1}.
\item 2003 - the non-symmetric version of the Finsleroid-Finsler metric in \cite{Vin9} as \emph{Asanov-type Finslerian metric functions}.     
\item 2006 - Asanov's examples for (non-symmetric) Finsleroid-Finsler metrics that are of Landsberg but not of Berwald type in \cite{Asanov2};
\item 2016 - non-symmetric Finsleroid-Finsler metrics with closed Finsleroid axis $1$-forms as the solutions of a conformal rigidity problem \cite{Vin18}.
\end{itemize}

In this paper we prove that Asanov's Unicorns belong to the class of generalized Berwald manifolds. More precisely we prove the following theorems: a Finsleroid-Finsler space is a generalized Berwald space if and only if the Finsleroid charge is constant. Especially a Finsleroid-Finsler space is a Landsberg space if and only if it is a generalized Berwald manifold with a semi-symmetric compatible linear connection. 

\section*{Acknowledgement}
The paper was motivated by the oral communication with Professor David Bao at the 50th Symposium on Finsler Geometry (21-25. Oct. 2015, Hiroshima, Japan). I am very grateful for his human and professional encouragement.

\section{Notations and terminology}

Let $M$ be a manifold with local coordinates $u^1, \ldots, u^n.$ The induced coordinate system of the tangent manifold $TM$ consists of the functions
$x^1, \ldots, x^n$ and $y^1, \ldots, y^n,$ where $x$'s refer to the coordinates of the base point and $y$'s denote the coordinates of the direction. A {\bf Finslerian metric}
is a continuous function $F\colon TM\to \mathbb{R}$ satisfying the following conditions:
$\displaystyle{F}$ is smooth on the complement of the zero section (regularity), $\displaystyle{F(tv)=tF(v)}$ for all $\displaystyle{t> 0}$ (positive homogenity) and the Hessian
$$g_{ij}=\frac{\partial^2 E}{\partial y^i \partial y^j}\ \ \textrm{of the Finslerian energy}\ \ E=\frac{1}{2}F^2$$
is positive definite at all nonzero elements $\displaystyle{v\in T_pM}$ (strong convexity). It is called the Riemann-Finsler metric of the Finsler manifold. The Riemann-Finsler metric makes each tangent space (except at the origin) a Riemannian manifold with standard canonical objects such as the {\bf volume form}
$$d\mu=\sqrt{\det g_{ij}}\ dy^1\wedge \ldots \wedge dy^n,$$
the {\bf Liouville vector field}  
$$C:=y^1\partial /\partial y^1 +\ldots +y^n\partial / \partial y^n$$
and the {\bf induced volume form on the indicatrix hypersurface} $\displaystyle{\partial K_p:=F^{-1}(1)\cap T_pM\ \  (p\in M)}$.
The coordinate expression is
$$\mu=\sqrt{\det g_{ij}}\ \sum_{i=1}^n (-1)^{i-1} \frac{y^i}{F} dy^1\wedge\ldots\wedge dy^{i-1}\wedge dy^{i+1}\ldots \wedge dy^n.$$
As a general reference of Finsler geometry see \cite{BCS}; we will use the following notations and terminology:
$$l_i=\frac{\partial F}{\partial y^i}, \ \ g^{ij}=(g_{ij})^{-1},\ \ C_{ijk}=\frac{1}{2}\frac{\partial g_{ij}}{\partial y^k},\ \ \textrm{where}\ \ \mathcal{C}^l_{ij}=g^{lk}\mathcal{C}_{ijk}$$
is the so-called {\bf first Cartan tensor}. The first Cartan tensor is totally symmetric and $\displaystyle{y^k\C_{ijk}=0}$. The {\bf geodesic spray coefficients} and the {\bf horizontal sections} are given by
\begin{equation}
\label{eq:3}
G^l=\frac{1}{2}g^{lm}\left(y^k\frac{\partial^2 E}{\partial y^m\partial x^k}-\frac{\partial E}{\partial x^m}\right)\ \ \textrm{and}\ \ \frac{\delta }{\delta x^i}=\frac{\partial }{\partial x^i}-G^l_i\frac{\partial }{\partial y^l},\ \ \textrm{where}\ \ G_i^l=\frac{\partial G^l}{\partial y^i}.
\end{equation}
The {\bf second Cartan tensor} or {\bf Landsberg tensor} is  
$$P^l_{ij}=\frac{1}{2} g^{lm}\left( \frac{\partial g_{jm}}{\partial x^i}-2G^k_i\C_{jkm}-G^k_{ij}g_{km}-G^k_{im} g_{jk}\right),\ \ \textrm{where}\ \ G_{ij}^l=\frac{\partial G_{i}^l}{\partial y^j}.$$
The {\bf mixed curvature of the Berwald connection} is defined as 
$$P_{ijk}^l=-G^{l}_{ijk}, \ \ \textrm{where}\ \ G_{ijk}^l=\frac{\partial G_{ij}^l}{\partial y^k}.$$
By some direct computations give the identity
\begin{equation}
\label{eq:4}
P^l_{ij}=-\frac{F}{2}l_m g^{kl}P_{ijk}^m.
\end{equation}

\begin{Def} Let $M$ be a Finsler manifold; a linear connection $\nabla$ on the base manifold $M$ is called \emph{compatible} to the Finslerian metric if the parallel transports with respect to $\nabla$ preserve the Finslerian lenght of tangent vectors. Finsler manifolds admitting compatible linear connections are called generalized Berwald manifolds. Berwald manifolds are generalized Berwald manifolds with torsion-free compatible linear connections. \end{Def}

\begin{Def} A Finsler manifold is a Landsberg manifold if the Landsberg tensor vanishes. 
\end{Def}

The notion of generalized Berwald manifolds goes back to V. Wagner \cite{Wag1}. The basic questions of the theory are the unicity of the compatible linear connection and its expression in terms of the canonical data of the Finsler manifold (intrinsic characterization). In case of a classical Berwald manifold (compatible linear connection with zero torsion) the intrinsic characterization is the vanishing of the mixed curvature tensor of the Berwald connection. This means that the quantities $G_{ij}^l$'s depend only on the position. They constitute the coefficients of the compatible linear connection on the base manifold. In general the intrinsic characterization of the compatible linear connection is based on the so-called {\bf averaged Riemannian metric}
\begin{equation}
\gamma_p(v,w):=\int_{\partial K_p}g(v,w)\, \mu.
\end{equation}

Using average processes is a new and important trend with a rapidly increasing number of papers in Finsler geometry; R. G. Torrom\'{e} \cite{Torrome1}, T. Aikou \cite{Aikou1}, M. Crampin \cite{C2} and \cite{C1},  V. S. Matveev (et. al.: H-B. Rademacher, M. Troyanov) \cite{Mat2}, \cite{Mat1} and \cite{Mat3}, Cs. Vincze \cite{Vin4}, \cite{Vin2}, \cite{Vin17} and \cite{Vin11}. For further references see also \cite{Vin13}, \cite{Vin14} and \cite{Vin12}.

\begin{Thm} \emph{\cite{Vin2}} If a linear connection on the base manifold is compatible with the Finslerian metric function then it must be metrical with respect to the averaged Riemannian metric $\gamma$.
\end{Thm}

It is well-known that a metric connection is uniquely determined by the torsion tensor. Following Agricola-Friedrich \cite{AF} consider the decomposition
$$T(X,Y):=T_1(X,Y)+T_2(X,Y),\ \ \textrm{where}\ \ \displaystyle{T_1(X,Y):=T(X,Y)-\frac{1}{n-1}\big(\tilde{T}(X)Y-\tilde{T}(Y)X\big)},$$
$\tilde{T}$ is the trace tensor of the torsion and
$$T_2(X,Y):=\frac{1}{n-1}\left(\tilde{T}(X)Y-\tilde{T}(Y)X\right).$$
Note that the torsion tensor of a metric linear connection on a manifold of dimension $2$ is automatically of the form $T=\tilde{T}(X)Y-\tilde{T}(Y)X$ (cf. Definition 3). Otherwise the trace-less part can be divided into two further components
$$T_1(X,Y)=A_1(X,Y)+S_1(X,Y)\ \ \Rightarrow\ \ T(X,Y)=A_1(X,Y)+S_1(X,Y)+T_2(X,Y)$$
by separating the totally anti-symmetric/axial part $A_1$. Therefore
we have eight possible classes of generalized Berwald manifolds depending on the surviving terms such as classical Berwald manifolds admitting torsion-free compatible linear connections \cite{Szab1} (we have no surviving terms) or Finsler manifolds admitting semi-symmetric compatible linear connections \cite{Vin17}.

\begin{Def} A linear connection is said to be semi-symmetric if the torsion tensor is of the form
\begin{equation}
\label{semi}
T(X,Y)=\lambda(Y)X-\lambda(X)Y
\end{equation}
where $\lambda$ is a one-form on the manifold. 
\end{Def}

The problem of the intrinsic characterization of compatible semi-symmetric linear connections is completely solved \cite{Vin17}: it can be expressed in terms of metrics and differential forms given by averaging. For the special case of exact (or at least closed) differential forms in the torsion (\ref{semi}) see \cite{Vin4}.

\begin{Thm} \emph{\cite{Vin17}}
A non-Riemannian Finsler manifold is a generalized Berwald manifold admitting a semi-symmetric compatible linear connection if and only if
$\sigma(p)>0$ for any $p\in M$ and the linear connection  
$$
{\bar{\nabla}}_X Y=\nabla^*_XY+\frac{1}{2\sigma}\left(\eta^*(Y)X-\gamma(X,Y)\eta^{* \sharp}\right)
$$
is compatible with the Finslerian metric function, where 
\begin{itemize}
\item $\nabla^*$ is the L\'{e}vi-Civita connection of the averaged Riemannian metric,
$$\ \ \ \ \rho^*:=\frac{d_{h^*}E}{E}-\frac{1}{2}\frac{S^*E}{E}\frac{d_J E^*}{E^*}\ \ \ \textrm{and}\ \ f:=\log \frac{\ E^*}{E},$$
where 
$$ J\left(\frac{\partial}{\partial x^i}\right)=\frac{\partial}{\partial y^i}\ \ \textrm{and}\ \ J\left(\frac{\partial}{\partial y^i}\right)=0$$ 
is the canonical vertical endomorphism/almost tangent structure on the tangent manifold,
\item $h^*$ is the associated horizontal endomorphism with $\nabla^*$; for the definition of the horizontal sections see 
formula \emph{(\ref{eq:3})} with substitution of the Riemannian energy $E^*$.
\item Furthermore
$$\eta^*(X_p):=\int_{\partial K^*_p} d_J\rho^*(\Theta, X^{h})+\frac{1}{2}\frac{S^* E}{E}X^vf\, \mu^*\ \ \ \textrm{and}\ \ \ \sigma(p):=\int_{\partial K^*_p} \frac{1}{2E^*}\|J\Theta\|^2\, \mu^*,$$
where 
$$\Theta=E^* g_*^{ij}\frac{\partial f}{\partial y^i} \frac{\partial}{\partial x^j}$$ 
is a gradient-type vector field and the norm is taken with respect to the vertically lifted Riemannian metric $g_{ij}^*=\gamma_{ij}\circ \pi$ and $g_*^{ij}=(g^*_{ij})^{-1}$.
\end{itemize}
\end{Thm}

\subsection{Randers metrics, ($\alpha, \beta$) - metrics, the sign property} The complete solution of the intrinsic characterization is also given in the special case of Randers manifolds without any special requirement for the torsion tensor. Let $(M, \alpha)$ be a connected Riemannian manifold and suppose that the one-form $\beta$ in $\wedge^1(M)$ satisfies condition
\begin{equation}
\sup_{\alpha(v,v)=1}\beta(v)< 1.
\end{equation}
The Randers metric on the manifold $M$ is defined  as $F(v)=\sqrt{\alpha(v,v)}+\beta(v).$

\begin{Thm} \emph{\cite{Vin1}} A Randers manifold is a generalized Berwald manifold if and only if there exists a linear connection $\nabla$ such that $\nabla \alpha=0$ and $\nabla \beta=0.$
\end{Thm}
The following theorem formulates a necessary and sufficient condition for a Randers manifold to be a generalized Berwald manifold in terms of the dual vector field
$\alpha(\beta^{\sharp},X)=\beta(X).$

\begin{Thm} \emph{\cite{Vin1}} A Randers manifold is a generalized Berwald manifold if and only if $\beta^{\sharp}$ is of constant Riemannian length.
\end{Thm}

The compatible linear connection is given as
\begin{equation}
\label{compatible}
\nabla_X Y=\nabla^*_X Y+\frac{\alpha(\nabla^*_X \beta^{\sharp}, Y)\beta^{\sharp}-\alpha(Y, \beta^{\sharp})\nabla^*_X {\beta^{\sharp}}}{\alpha(\beta^{\sharp}, \beta^{\sharp})}
\end{equation}

\noindent
If the compatible linear connection is semi-symmetric then we also have a structure theorem for the Riemannian manifold admitting a perturbation $\beta$ such that the Randers manifold is a generalized Berwald manifold with a semi-symmetric compatible linear connection \cite{Vin1}. By the main result in \cite{Vin1} the manifold carries a warped product metric structure; for the special case of exact (or at least closed) differential forms in the torsion (\ref{semi}) see \cite{Vin16}, see also \cite{Vin7}. These results have been generalized by Tayebi and Barzegari in \cite{Tayebi} for $(\alpha, \beta)$ -   metrics satisfying the {\bf{sign property}}
\begin{equation}
\label{signproperty}
\Phi'(-s)\Phi(s)+\Phi(-s)\Phi'(s) > 0 \ \ \textrm{or}\ \ \Phi'(-s)\Phi(s)+\Phi(-s)\Phi'(s) < 0,
\end{equation}
where $F=\alpha\Phi\left(\frac{\beta}{\alpha}\right)$ is a Finslerian metric function and $\Phi\colon (-b_0, b_0)\to \mathbb{R}^+.$ 
According to the positivity of $\Phi$ the sign property (\ref{signproperty}) is equivalent to
\begin{equation}
\label{signproperty1}
\varphi'(-s)\varphi(s)+\varphi(-s)\varphi'(s)> 0 \ \ \textrm{or}\ \ \varphi'(-s)\varphi(s)+\varphi(-s)\varphi'(s)> 0\ \ \textrm{where}\ \ \varphi=\Phi^2.
\end{equation}

\begin{Thm} \emph{\cite{Tayebi}} A Finsler manifold with an $(\alpha, \beta)$ - metric satisfying the sign property is a generalized Berwald manifold if and only if there exists a linear connection $\nabla$ such that $\nabla \alpha=0$ and $\nabla \beta=0.$
\end{Thm}

\begin{Thm} \emph{\cite{Tayebi}} A Finsler manifold with an $(\alpha, \beta)$ - metric satisfying the sign property is a generalized Berwald manifold if and only if $\beta^{\sharp}$ is of constant Riemannian length.
\end{Thm}

\section{Asanov's Finsleroid-Finsler metrics}
Using Asanov's original notations in \cite{Asanov2} the general form of Finsleroid-Finsler metrics is given by
\begin{equation}
\label{asanovmetric1}
F=e^{G\Phi/2}\sqrt{b^2+gqb+q^2},
\end{equation}
where $b=b_iy^i$ is the Finsleroid axis $1$ - form, $q=\sqrt{r_{ij}y^iy^j}$, $r_{ij}=a_{ij}-b_ib_j$ and $a_{ij}$ is a Riemannian metric such that $a^{ij}b_ib_j=1$, 
$$g=g(p)\ \ \textrm{and}\ \ -2 < g(p)< 2 \ \ (\textrm{the Finsleroid charge}),$$
$$h=\sqrt{1-\frac{\ g^2}{4}}, \ \ G=g/h,\ \ \Phi=\left\{\begin{array}{rl}
&+\frac{\pi}{2}+\arctan \frac{G}{2}-\arctan \frac{q+\frac{g}{2}b}{hb}\ \ \textrm{if} \ \ b > 0\\
&\\
&-\frac{\pi}{2}+\arctan \frac{G}{2}-\arctan \frac{q+\frac{g}{2}b}{hb}\ \ \textrm{if} \ \ b < 0.
\end{array}
\right.$$
The common limit of the right hand sides as $b\to 0$ is $\displaystyle{\arctan \frac{G}{2}}$. 

\subsection{An alternative formulation} In what follows we present the metric in a more compact form. If
$$f(b):=\arctan \frac{q+\frac{g}{2}b}{hb}$$
then
\begin{equation}
\label{f}
\lim_{b\to \pm \infty}f(b)=\arctan \frac{g}{2h}=\arctan \frac{G}{2}\ \ \textrm{and}\ \ f'(b)=-\frac{qh}{b^2+q^2+bgq}
\end{equation}
because $h^2+g^2/4=1$. In a similar way, if 
$$\tilde{f}(b):=\arctan \frac{2b+gq}{2hq}$$ 
then
\begin{equation}
\label{h}
\lim_{b\to \pm \infty}\tilde{f}(b)=\pm \frac{\pi}{2}\ \ \textrm{and}\ \ \tilde{f}'(b)=\frac{qh}{b^2+q^2+bgq}.
\end{equation}
Therefore $f+\tilde{f}$ is constant on the connected parts of the domain. Taking the limits $b\to \infty$ and $b\to -\infty$, respectively, we have 
$$\arctan \frac{q+\frac{g}{2}b}{hb}+\arctan \frac{2b+gq}{2hq}=\left\{\begin{array}{rl}
&+\frac{\pi}{2}+\arctan \frac{G}{2}\ \ \textrm{if} \ \ b>0\\
&\\
&-\frac{\pi}{2}+\arctan \frac{G}{2}\ \ \textrm{if} \ \ b < 0.
\end{array}
\right.\ \ \Rightarrow$$
\begin{equation}
\label{metrics}
\Phi=\arctan \frac{2b+gq}{2hq}.
\end{equation}

\begin{Def}
Using the notations
$$\alpha=\sqrt{a_{ij}y^iy^j},\ \ \beta=b_iy^i\ \ (\textrm{Finsleroid axis 1 - form})\ \ \textrm{and}\ \ g=\frac{K}{2} \ \ (\textrm{Finsleroid charge})$$
it follows that  
$$q=\sqrt{\alpha^2-\beta^2},\ \ h:=\sqrt{1-\frac{\ g^2}{4}}=\sqrt{1-\frac{\ K^2}{16}}=\frac{\sqrt{16-K^2}}{4}\ \ \textrm{and}\ \ G=g/h=\frac{2K}{\sqrt{16-K^2}}$$
and the Finslerian energy $E=(1/2)F^2$ of a Finsleroid-Finsler metric is
\begin{equation}
\label{asanovmetric}
E=\frac{1}{2}\left(\alpha^2+\frac{K}{2}\beta\sqrt{\alpha^2-\beta^2}\right)e^{\frac{2K}{\sqrt{16-K^2}}\arctan \frac{1}{\sqrt{16-K^2}}\left(\frac{4\beta}{\sqrt{\alpha^2-\beta^2}}+K\right)}.
\end{equation}
\end{Def}

The metric (\ref{asanovmetric}) is formally an $(\alpha, \beta)$ - metric
\begin{equation}
\label{alphabeta}
E=\frac{1}{2}\alpha^2 \varphi \left(\frac{\beta}{\alpha}\right),
\end{equation}
where 
\begin{equation}
\label{phi1}
\varphi(s)=\left(1+\frac{K}{2}s\sqrt{1-s^2}\right)e^{\frac{2K}{\sqrt{16-K^2}} \arctan \frac{1}{\sqrt{16-K^2}}\left(\frac{4s}{\sqrt{1-s^2}}+K\right)},
\end{equation}
$s\in [-1,1]$ and the value at $\pm 1$ is defined by the continuous extension
\begin{equation}
\label{cont}
\varphi(\pm 1)=\lim_{s\to \pm 1}\varphi(s)= e^{\pm\frac{2K}{\sqrt{16-K^2}}\frac{\pi}{2}}\ \ \Rightarrow\ \ E(\pm \beta^{\sharp})= \frac{1}{2}e^{\pm\frac{2K}{\sqrt{16-K^2}}\frac{\pi}{2}}.
\end{equation}
Actually (\ref{alphabeta}) represents a more general form of metrics because $\varphi$ depends on the position too. In case of a standard ($\alpha, \beta$) - metric $\varphi$ is a  function of the single variable $s$.

\begin{Lem} The function $\varphi$ is of class $\mathcal{C}^1$ with respect to the variable $s$.
\end{Lem}

\begin{Pf}
Fix a point $p\in M$; in what follows we prove that $\varphi$ is of class $\mathcal{C}^1$ at $\pm 1$. We discuss the case of $s=1$ in details. The case of $s=-1$ is similar. By definition
$$\varphi'(1)=\lim_{s\to 1}\frac{\varphi(s)-\varphi(1)}{s-1}.$$
For any fixed $s<1$ we can use the Lagrange mean value theorem 
$$\varphi(s)-\varphi(1)=\varphi'(t)(s-1)$$
because $\varphi$ is continuous on the closed inteval $[s,1]$ (see formula (\ref{cont}) of the continuous extension) and differentiable on $]s,1[$. Therefore
\begin{equation}
\label{form1}
\varphi'(1)=\lim_{t\to 1}\varphi'(t)=0
\end{equation}
as a simple calculation shows (see Figure 1):
\begin{equation}
\label{phi2}
\varphi'(s)=K\sqrt{1-s^2}e^{\frac{2K}{\sqrt{16-K^2}} \arctan \frac{1}{\sqrt{16-K^2}}\left(\frac{4s}{\sqrt{1-s^2}}+K\right)}
\end{equation}
\end{Pf}

\begin{figure}
\centering
\includegraphics[viewport=0 0 474 172, scale=0.8]{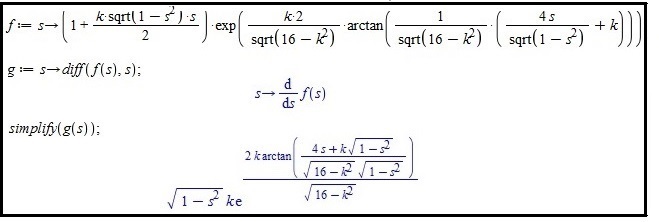}
\caption{The first derivative.}
\end{figure}

\begin{Thm} \label{tsign}
The function $\varphi$ satisfies the sign property 
$$a(s):=\varphi'(-s)\varphi(s)+\varphi(-s)\varphi(-s)=2K\sqrt{1-s^2}e^{2A(-s)+2A(s)} \left\{\begin{array}{rl}
& > 0\ \ \textrm{if}\ \ K > 0,\\
& \\
& < 0\ \ \textrm{if}\ \ K < 0,
\end{array}
\right.$$
where
$$2A(s):=\frac{2K}{\sqrt{16-K^2}} \arctan \frac{1}{\sqrt{16-K^2}}\left(\frac{4s}{\sqrt{1-s^2}}+K\right).$$
\end{Thm}

\begin{Pf}
The proof is a straightforward calculation as Figure 2 shows (the worksheet is the continuation of Figure 1).
\end{Pf}

\begin{figure}
\centering
\includegraphics[viewport=0 0 448 262, scale=0.8]{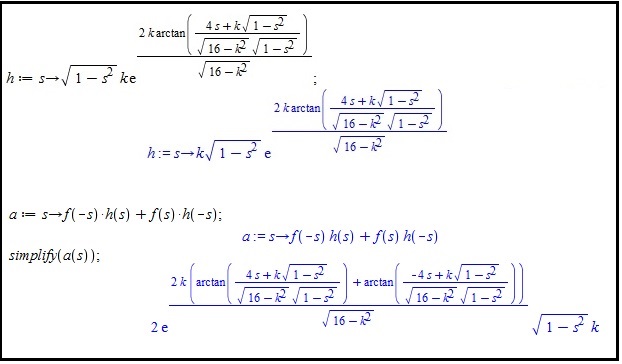}
\caption{The sign-property.}
\end{figure}

\begin{Thm} \emph{\cite{Asanov2}} \label{ltheorem} A Finsleroid-Finsler space is a Landsberg space if and only if the function $K$ is constant and
\begin{equation}
\label{landsberg}
(\nabla^* \beta)(X,Y)=\frac{\textrm{\ \emph{div}}\ \beta^{\sharp} }{n-1}\left(\alpha(X,Y)-\beta(X)\beta(Y)\right)
\end{equation}
\end{Thm}

\begin{Rem}
{\emph{The original formulation of Theorem \ref{ltheorem} in \cite{Asanov2} (Theorem 3, p. 278) is that a Finsleroid-Finsler space is a Landsberg space if and only if the Finsleroid axis 1-form $\beta$ is closed, the Finsleroid charge is constant and
\begin{equation}
\label{landsberg1}
(\nabla^* \beta)(X,Y)=k\left(\alpha(X,Y)-\beta(X)\beta(Y)\right)
\end{equation}
for some scalar field $k\colon M\to \mathbb{R}$. Note that if (\ref{landsberg1}) holds then the closedness of $\beta$ is redundant because
$$d\beta (X,Y)=\nabla^*\beta (X,Y)-\nabla^*\beta (Y,X)\stackrel{(\ref{landsberg1})}{=}0.$$
On the other hand condition (\ref{landsberg1}) is obviously equivalent to
$$\alpha(\nabla^*_X\beta^{\sharp}, Y)=k\left(\alpha(X,Y)-\beta(X)\beta(Y)\right)\ \ \Rightarrow\ \ k=\frac{\ \textrm{div}\ \beta^{\sharp}}{n-1}$$
because of the unit length of $\beta^{\sharp}$ with respect to $\alpha$.}}
\end{Rem}

\section{The main results}

\begin{Thm}
\label{t1}
A connected Finsleroid-Finsler space is a generalized Berwald space if and only if the Finsleroid charge is constant.
\end{Thm}

\begin{Pf}
Suppose that there exists a linear connection such that it is compatible to the Finsleroid-Finsler metric. Since the parallel transports preserve the Finslerian norm of tangent vectors and they are linear between the tangent spaces it follows that they preserve the Riemann-Finsler metric $g_{ij}$ and the indicatrices with the induced Riemannian metric are isometric. Asanov's Finsleroid-Finsler metric has indicatrices of constant curvature 
$$1-\frac{K^2(p)}{4};$$
see \cite{Asanov2}, formula (2.32). Therefore the Finsleroid charge $K$ must be constant on a connected manifold. Conversely, suppose that the function $K$ is constant. Then the Finsleroid-Finsler metric is an $(\alpha, \beta)$ - metric of the form
$$F=\alpha \phi\left(\frac{\beta}{\alpha}\right),\ \ \textrm{where}\ \ \phi^2=\varphi$$
because $\varphi$ does not depend on the position; see formula (\ref{alphabeta}) and Definition 4. We have two possible cases: 
\begin{itemize}
\item If $K=0$ then the space is Riemannian as a special generalized Berwald manifold.
\item If $K\neq 0$ then Theorem \ref{tsign} implies that the function $\varphi$ satisfies the sign property (\ref{signproperty1}) and, by Theorem 6, it is a generalized Berwald manifold because the Finsleroid axis 1-form $\beta$ is of constant (unit) length with respect to $\alpha$.
\end{itemize}
\end{Pf}

\begin{Thm}
A connected Finsleroid-Finsler space is a Landsberg space if and only if it is a generalized Berwald space with a semi-symmetric compatible linear connection such that the torsion tensor is of the form 
$$T(X,Y)=\lambda(Y)X-\lambda(X)Y,\ \ \textrm{where}\ \ \lambda=-\frac{\textrm{\emph{div}}\ \beta^{\sharp}}{n-1}\beta.$$
\end{Thm}

\begin{Pf}
Suppose that a Finsleroid-Finsler space is a Landsberg space. Then, by Theorem 8, we have that the Finsleroid charge $K$ is constant, i.e. the space is a generalized Berwald space in the sense of Theorem 9. By formula (\ref{compatible}) the compatible linear connection is
$$\nabla_X Y=\nabla^*_X Y+\frac{\alpha(\nabla^*_X \beta^{\sharp}, Y)\beta^{\sharp}-\alpha(Y, \beta^{\sharp})\nabla^*_X {\beta^{\sharp}}}{\alpha(\beta^{\sharp}, \beta^{\sharp})}=$$
$$\nabla^*_X Y+\alpha(\nabla^*_X \beta^{\sharp}, Y)\beta^{\sharp}-\alpha(Y, \beta^{\sharp})\nabla^*_X {\beta^{\sharp}}$$
because $\beta^{\sharp}$ is of unit length with respect to $\alpha$. On the other hand
$$\alpha(\nabla^*_X \beta^{\sharp}, Y)=(\nabla^* \beta) (X, Y)\stackrel{(\ref{landsberg})}{=}\frac{\textrm{\ div}\ \beta^{\sharp} }{n-1}\left(\alpha(X,Y)-\beta(X)\beta(Y)\right)\ \ \Rightarrow\ \ $$
$$\nabla^*_X \beta^{\sharp}=\frac{\ \textrm{div}\  \beta^{\sharp} }{n-1}\left(X-\beta(X)\beta^{\sharp}\right)$$
and, consequently,
\begin{equation}
\label{formula1}
\nabla_X Y=\nabla^*_X Y-\frac{\ \textrm{div}\  \beta^{\sharp}}{n-1}\beta(Y)X+\alpha(X,Y)\frac{\ \textrm{div}\  \beta^{\sharp}}{n+1} \beta^{\sharp},
\end{equation}
i.e.
\begin{equation}
\label{formula2}
\nabla_X Y=\nabla^*_X Y+\lambda(Y)X-\alpha(X,Y)\lambda^{\sharp},\ \ \textrm{where}\ \ \lambda=-\frac{\ \textrm{div}\  \beta^{\sharp}}{n-1}\beta. 
\end{equation}
Formula (\ref{formula2}) determines the only metric linear connection with torsion 
$$T(X,Y)=\lambda(Y)X-\lambda(X)Y.$$
Conversely, suppose that we have a Finsleroid-Finsler space such that it is a generalized Berwald manifold with $\nabla$ in formula (\ref{formula1}) as a compatible linear connection. Then, by Theorem \ref{t1}, the Finsleroid charge is constant and we have an ($\alpha, \beta$) - metric. We have two possible cases: 
\begin{itemize}
\item If $K=0$ then the space is Riemannian as a special Landsberg manifold.
\item If $K\neq 0$ then Theorem \ref{tsign} implies that the function $\varphi$ satisfies the sign property (\ref{signproperty1}) and, by Theorem 5, $\nabla \beta=0$. This implies the special form (\ref{landsberg}) of $\nabla^* \beta$ and the statement follows by Theorem \ref{ltheorem}.
\end{itemize} 
\end{Pf}

\section{Appendix: regularity properties of Finsleroid-Finsler metrics}

Finsleroid-Finsler metrics belong to the class of $y$ - local Finslerian metrics because the third order partial derivatives with respect to the variables $y$'s are singular at $\pm \beta^{\sharp}$. In what follows we prove that the partial derivatives with respect to the variables $y$'s exist and continuous up to order $2$, i.e. Finsleroid-Finsler metrics are of class $C^2$ on the complement of the zero section. 

\subsection{The first $y$ - derivatives of a Finsleroid-Finsler energy function} By Lemma 1 the function $\varphi$ is of class $C^1$. Therefore
\begin{equation}
\label{first1}
\frac{\partial E}{\partial y^i}\stackrel{(\ref{alphabeta})}{=}\frac{1}{2}\left(2y^i\varphi\left(\frac{\beta}{\alpha}\right)+\alpha^2\varphi'\left(\frac{\beta}{\alpha}\right)\frac{\partial \beta/\alpha}{\partial y^i}\right)=\frac{1}{2}\left(2y^i\varphi\left(\frac{\beta}{\alpha}\right)+\alpha^2\varphi'\left(\frac{\beta}{\alpha}\right)\frac{b_i\alpha-(\beta/\alpha)y^i}{\alpha^2}\right)
\end{equation}
is of class $\mathcal{C}^1$ on the complement of the zero section. For the sake of simplicity we can use an orthonormal coordinate system 
\begin{equation}
\label{orto}
y^1, \ldots, y^n\ \ \textrm{such that}\ \ y^n=\beta \ \ \textrm{i.e.}\ \ b_i=\delta_i^n \ \ \textrm{and}\ \ \alpha^2=(y^1)^2+\ldots+(y^n)^2.
\end{equation}
In terms of an orthonormal coordinate system (\ref{orto})
$$\frac{\partial E}{\partial y^i}\stackrel{(\ref{phi1}), (\ref{phi2}), (\ref{first1})}{=}e^{\frac{2K}{\sqrt{16-K^2}} \arctan \frac{1}{\sqrt{16-K^2}}\left(\frac{4\left(\frac{\beta}{\alpha}\right)}{\sqrt{1-\left(\frac{\beta}{\alpha}\right)^2}}+K\right)}\left(y^i+\frac{K}{2}\delta_i^n \alpha \sqrt{1-\left(\frac{\beta}{\alpha}\right)^2}\right)\ \ \Rightarrow$$
\begin{equation}
\label{first2}
\frac{\partial E}{\partial y^i}=e^{\frac{2K}{\sqrt{16-K^2}} \arctan \frac{1}{\sqrt{16-K^2}}\left(\frac{4\beta}{\sqrt{\alpha^2-\beta^2}}+K\right)}\left(y^i+\frac{K}{2}\delta_i^n \sqrt{\alpha^2-\beta^2}\right).
\end{equation}
By (\ref{form1}) and (\ref{first1})
\begin{equation}
\label{first3}
\frac{\partial E}{\partial y^i}_{(0, \ldots, 0, \pm1)}=\delta_i^n e^{\pm\frac{2K}{\sqrt{16-K^2}} \frac{\pi}{2}}.
\end{equation}

\subsection{The second $y$ - derivatives of a Finsleroid-Finsler energy function}
To compute the second order partial derivatives it is useful to introduce the function
$$\ \ 2A:=\frac{2K}{\sqrt{16-K^2}} \arctan \frac{1}{\sqrt{16-K^2}}\left(\frac{4\beta}{\sqrt{\alpha^2-\beta^2}}+K\right).$$
It can be expressed as the function of the variable $s=\beta/\alpha$
$$2A(s):=\frac{2K}{\sqrt{16-K^2}} \arctan \frac{1}{\sqrt{16-K^2}}\left(\frac{4s}{\sqrt{1-s^2}}+K\right).$$

\noindent
Since formula (\ref{first2}) can be written as 
$$\frac{\partial E}{\partial y^i}=e^{2A\left(\frac{\beta}{\alpha}\right)}\left(y^i+\frac{K}{2}\delta_i^n\sqrt{\alpha^2-\beta^2}\right) $$
\begin{figure}
\centering
\includegraphics[viewport=0 0 365 228, scale=0.8]{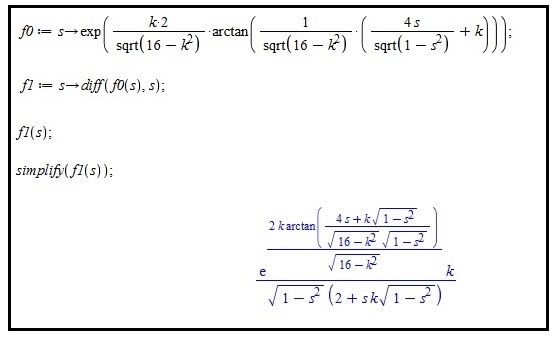}
\caption{The derivative of the exponential term}
\end{figure}
it follows that 
\begin{equation}
\label{second0}
\frac{\partial^2 E}{\partial y^j \partial y^i}=\frac{\partial e^{2A}}{\partial y^j}\left(y^i+\frac{K}{2}\delta_i^n\sqrt{\alpha^2-\beta^2}\right)+e^{2A\left(\frac{\beta}{\alpha}\right)}\left(\delta_j^i+\frac{K}{2}\delta_i^n\frac{y^j-\beta \delta_j^n}{\sqrt{\alpha^2-\beta^2}}\right),
\end{equation}
where
$$\frac{\partial e^{2A}}{\partial y^j}=e^{2A\left(\frac{\beta}{\alpha}\right)}(2A)'\left(\frac{\beta}{\alpha}\right)\frac{\partial \beta/\alpha}{\partial y^j}=K\frac{e^{2A\left(\frac{\beta}{\alpha}\right)}}{2\sqrt{1-\left(\frac{\beta}{\alpha}\right)^2}\left(1+\frac{K}{2}\frac{\beta}{\alpha}\sqrt{1-\left(\frac{\beta}{\alpha}\right)^2}\right)}\frac{b_j\alpha-(\beta/\alpha)y^j}{\alpha^2}$$
as a straightforward calculation shows; see Figure 3. Therefore
\begin{equation}
\label{expder}
\frac{\partial e^{2A}}{\partial y^j}=K\frac{e^{2A\left(\frac{\beta}{\alpha}\right)}}{2\sqrt{\alpha^2-\beta^2}\left(\alpha^2+\frac{K}{2}\beta\sqrt{\alpha^2-\beta^2}\right)}(\delta^n_j\alpha^2-\beta y^j)
\end{equation}
and we have
$$
\frac{\partial^2 E}{\partial y^j \partial y^i}\stackrel{(\ref{second0}), (\ref{expder})}{=}\frac{e^{2A\left(\frac{\beta}{\alpha}\right)}}{\alpha^2+\frac{K}{2}\beta\sqrt{\alpha^2-\beta^2}}
$$
$$\left( K\frac{\alpha^2\delta_j^n-\beta y^j}{2\sqrt{\alpha^2-\beta^2}}\left(y^i+\frac{K}{2}\delta_i^n\sqrt{\alpha^2-\beta^2}\right)+\left(\delta_j^i+\frac{K}{2}\delta_i^n\frac{y^j-\beta\delta_j^n}{\sqrt{\alpha^2-\beta^2}}\right)\left(\alpha^2+\frac{K}{2}\beta\sqrt{\alpha^2-\beta^2}\right)\right),$$
where
{\small{$$ K\frac{\alpha^2\delta_j^n-\beta y^j}{2\sqrt{\alpha^2-\beta^2}}\left(y^i+\frac{K}{2}\delta_i^n\sqrt{\alpha^2-\beta^2}\right)+\left(\delta_j^i+\frac{K}{2}\delta_i^n\frac{y^j-\beta\delta_j^n}{\sqrt{\alpha^2-\beta^2}}\right)\left(\alpha^2+\frac{K}{2}\beta\sqrt{\alpha^2-\beta^2}\right)=$$
$$K\frac{(\alpha^2\delta_j^n-\beta y^j)y^i+\delta_i^n(y^j-\beta\delta_j^n)\alpha^2}{2\sqrt{\alpha^2-\beta^2}}+\delta_i^j\left(\alpha^2+\frac{K}{2}\beta\sqrt{\alpha^2-\beta^2}\right)+\frac{K^2}{4}(\delta_j^n\alpha^2-\beta y^j)\delta_i^n+\frac{K^2}{4}\beta \delta_i^n(y^j-\beta\delta_j^n)=$$
$$K\frac{\alpha^2(\delta_j^ny^i+\delta_i^ny^j)-\beta(y^iy^j+\alpha^2\delta_i^n \delta_j^n)}{2\sqrt{\alpha^2-\beta^2}}+\delta_i^j\left(\alpha^2+\frac{K}{2}\beta\sqrt{\alpha^2-\beta^2}\right)+\frac{K^2}{4}\delta_i^n \delta_j^n(\alpha^2-\beta^2).$$}}

\noindent
Finally
\begin{equation}
\label{second}
\frac{\partial^2 E}{\partial y^j \partial y^i}=\frac{e^{2A\left(\frac{\beta}{\alpha}\right)}}{\alpha^2+\frac{K}{2}\beta\sqrt{\alpha^2-\beta^2}}
\end{equation}
$$\left( K\frac{\alpha^2(\delta_j^ny^i+\delta_i^ny^j)-\beta(y^iy^j+\alpha^2\delta_i^n \delta_j^n)}{2\sqrt{\alpha^2-\beta^2}}+\delta_i^j\left(\alpha^2+\frac{K}{2}\beta\sqrt{\alpha^2-\beta^2}\right)+\frac{K^2}{4}\delta_i^n \delta_j^n(\alpha^2-\beta^2)\right).$$

\subsection{The continuity of the second order partial derivatives at $(0, \ldots, 0, \pm 1)$}

\subsubsection{The case of $i=1, \ldots, n-1$ and $j=1, \ldots, n-1$}

If $1 \leq i, j < n$ then by (\ref{second})
\begin{equation}
\label{second1}
\frac{\partial^2 E}{\partial y^j \partial y^i}=\frac{e^{2A\left(\frac{\beta}{\alpha}\right)}}{\alpha^2+\frac{K}{2}\beta\sqrt{\alpha^2-\beta^2}}\left(\frac{-K\beta y^iy^j}{2\sqrt{\alpha^2-\beta^2}}+\delta_i^j\left(\alpha^2+\frac{K}{2}\beta\sqrt{\alpha^2-\beta^2}\right)\right).
\end{equation}
On the other hand 
$$\left|\frac{y^i}{\sqrt{\alpha^2-\beta^2}}\right|=\left|\frac{y^i}{\sqrt{(y^1)^2+\ldots+(y^i)^2+\ldots (y^{n-1})^2}}\right| \leq 1 \ \ (\textrm{boundedness})$$
and, consequently,
$$\frac{-K\beta y^iy^j}{2\sqrt{\alpha^2-\beta^2}}\ \mapsto \ \ 0\ \ \textrm{as} \ \ y^j \to 0.$$
Therefore 
\begin{equation}
\label{limit}
\lim_{v\to (0, \ldots, 0, \pm 1)}\frac{\partial^2 E}{\partial y^j \partial y^i}_v\stackrel{(\ref{second1})}{=}e^{\pm \frac{2K}{\sqrt{16-K^2}}\frac{\pi}{2}}\delta_i^j.
\end{equation}
By definition
$$\frac{\partial^2 E}{\partial y^j \partial y^i}_{(0, \ldots, 0, \pm 1)}:=\lim_{s\to 0}\frac{\frac{\partial E}{\partial y^i}_{(0, \ldots, 0, s, 0, \ldots, 0, \pm 1)}-\frac{\partial E}{\partial y^i}_{(0, \ldots, 0,\pm 1)}}{s}.$$
Using that the first order partial derivatives are continuous (subsection 4.1) the Lagrange mean value theorem shows that the second order partial derivatives at $(0, \ldots, 0, \pm 1)$ is just the limit (\ref{limit}):
$$\frac{\partial E}{\partial y^i}_{(0, \ldots, 0, s, 0, \ldots, 0, \pm 1)}-\frac{\partial E}{\partial y^i}_{(0, \ldots, 0, \pm 1)}=\frac{\partial^2 E}{\partial y^j \partial y^i}_{(0, \ldots, 0, t, 0, \ldots, 0, \pm 1)}(s-0),$$
where $t$ is between $0$ and $s$, i.e.
$$\frac{\partial^2 E}{\partial y^j \partial y^i}_{(0, \ldots, 0,\pm 1)}=\lim_{t\to 0}\frac{\partial^2 E}{\partial y^j \partial y^i}_{(0, \ldots, 0, t, 0, \ldots, 0, \pm 1)}\stackrel{(\ref{limit})}{=}e^{\pm \frac{2K}{\sqrt{16-K^2}}\frac{\pi}{2}}\delta_i^j$$
and we have the continuity of the second order partial derivatives.

\subsubsection{The case of $i=1, \ldots, n-1$ and $j=n$} Using that $y^n=\beta$ formula (\ref{second}) implies that 
$$\frac{\partial^2 E}{\partial y^n \partial y^i}=\frac{e^{2A\left(\frac{\beta}{\alpha}\right)}}{\alpha^2+\frac{K}{2}\beta \sqrt{\alpha^2-\beta^2}}\left(\frac{K\alpha^2 y^i-K\beta^2 y^i}{2\sqrt{\alpha^2-\beta^2}}+\delta_i^n \left( \alpha^2+\frac{K}{2}\beta\sqrt{\alpha^2-\beta^2}\right)\right)=$$
$$\frac{e^{2A}}{\alpha^2+\frac{K}{2}\beta\sqrt{\alpha^2-\beta^2}}\frac{K}{2}y^i\sqrt{\alpha^2-\beta^2}.$$
Therefore
$$\lim_{v\to (0, \ldots, 0,\pm 1)}\frac{\partial^2 E}{\partial y^n \partial y^i}_v=0.$$
The computation of the second order partial derivative at $(0, \ldots, 0, \pm 1)$ by definition needs the same step based on the Lagrange mean value theorem as in subsection 4.3.1.

\subsubsection{The case of $i=j=n$} Since $y^n=\beta$ it follows by (\ref{second}) that
$$\frac{\partial^2 E}{\partial y^n \partial y^n}=\frac{e^{2A\left(\frac{\beta}{\alpha}\right)}}{\alpha^2+\frac{K}{2}\beta\sqrt{\alpha^2-\beta^2}}\left(\frac{2K\alpha^2 \beta -K\beta(\beta^2+\alpha^2)}{2\sqrt{\alpha^2-\beta^2}}+\alpha^2+\frac{K}{2}\sqrt{\alpha^2-\beta^2}\beta +\frac{K^2}{4}(\alpha^2-\beta^2)\right)=$$
$$\frac{e^{2A}}{\alpha^2+\frac{K}{2}\beta\sqrt{\alpha^2-\beta^2}}\left(K\beta\sqrt{\alpha^2-\beta^2}+\alpha^2 +\frac{K^2}{4}(\alpha^2-\beta^2)\right).$$
Therefore
$$\lim_{v\to (0, \ldots, 0,\pm 1)}\frac{\partial^2 E}{\partial y^n \partial y^n}_v=e^{\pm \frac{2K}{\sqrt{16-K^2}}\frac{\pi}{2}}.$$
The computation of the second order partial derivative at $(0, \ldots, 0, \pm 1)$ by definition needs the same step based on the Lagrange mean value theorem as in subsection 4.3.1.

\end{document}